\documentclass[12pt]{amsart}
\allowdisplaybreaks
\usepackage{amssymb}
\usepackage{graphicx}
\usepackage{fullpage}
\usepackage[pdfencoding=auto,colorlinks,linkcolor=,citecolor=]{hyperref}
\usepackage{cite}
\usepackage{dsfont}
\usepackage{eucal}
\usepackage[all,cmtip]{xy}
\newcommand{\hra}{\hookrightarrow}

\newcommand{\phat}{{}^{^\wedge}_p}
\newcommand{\Qhat}{{}^{^\wedge}_\bQ}

\newcommand{\bC}{\mathbb C}
\newcommand{\bF}{\mathbb F}
\newcommand{\bi}{{\mathsf{i}\mkern1mu}}
\newcommand{\bQ}{\mathbb Q}

\newcommand{\bZ}{\mathbb Z}

\newcommand{\kk}{\mathsf{k}}

\renewcommand{\ge}{\geqslant}

\makeindex
\numberwithin{equation}{section}  

\theoremstyle{plain}
\newtheorem{lemma}[equation]{Lemma}
\newtheorem{theorem}[equation]{Theorem}
\newtheorem{proposition}[equation]{Proposition}

\theoremstyle{definition}

\theoremstyle{remark} 
\newtheorem{remark}[equation]{Remark}

\author{David J. Benson} 
\address{Institute of Mathematics \\ 
University of Aberdeen \\ 
Aberdeen AB24 3UE \\ 
United Kingdom}
\email{d.j.benson@abdn.ac.uk}

\author{John Greenlees} 
\address{Mathematics Institute, Zeeman Building, University of
  Warwick, Coventry CV4 7AL, United Kingdom}
\email{John.Greenlees@warwick.ac.uk}

\title{Formality of cochains on $BG$}

\begin{document}

\begin{abstract}
Let $G$ be a compact Lie group with maximal torus $T$.
If $|N_G(T)/T|$ is invertible in the field $\kk$ then
the algebra of cochains $C^*(BG;\kk)$ is formal as an $A_\infty$ algebra,
or equivalently as a DG algebra.
\end{abstract}

\maketitle

\section{Introduction}

It is well known that over a field of characteristic zero, the algebra
of cochains 
on the classifying space of a connected compact Lie group is formal
as an $A_\infty$ algebra, or equivalently as a differential graded
(DG) algebra. 
We prove that this is the case for a compact Lie group $G$ that is not
necessarily connected, over any field in which the order of $N_G(T)/T$
 is invertible, where $T$ is a maximal torus in $G$. 

In contrast, 
even for a finite group with cyclic Sylow subgroups 
of order $p^n\ge 3$ in characteristic $p$, the algebra of
 cochains on the classifying space is not formal. The $A_\infty$
 structure in this case is computed in~\cite{Benson/Greenlees:2021a}.

Our notation is as follows. We write $G$ for a compact Lie group,
$T$ for a maximal torus in $T$, and $W$ for the finite group 
$N_G(T)/T$. This a finite group acting on $T$ by conjugation, and
the kernel of the action is $C_G(T)/T$. If $G$ is connected then
$C_G(T)=T$, and $W$ is the Weyl group of $G$.
If $p$ is a prime, we write $X\phat$ for the
Bousfield--Kan $p$-completion of a space $X$,
see~\cite{Bousfield/Kan:1972a}. 

Our strategy is quite different to the usual characteristic zero proof. 
After reducing to the case of a normal torus, 
we make an approximation to $G$ at the prime $p$ 
by a locally finite group $\breve G$ with $B\breve G\phat$ homotopy
equivalent to $BG\phat$, and then we put an internal grading on the 
group algebra $\bF_p\breve G$. This gives us a second, internal grading on 
cochains  $C^*(B\breve G;\bF_p)$, and the $A_\infty$ structure maps $m_i$
provided by Kadeishvili's theorem on $H^*(B\breve G;\bF_p)$ have to preserve
the internal grading. This then proves that they are all zero apart from 
the multiplication map $m_2$. Our main theorem is the following.

\begin{theorem}\label{th:main}
If $G$ is a compact Lie group with maximal torus $T$,  and $p$ does
not divide $|N_G(T)/T|$, then
$C^*(BG;\bF_p)$ is a formal $A_\infty$ algebra.
\end{theorem}

In Section~\ref{se:0} we recall the proof of the characteristic zero
theorem. The proof of Theorem~\ref{th:main} occupies
Sections~\ref{se:finite} and~\ref{se:formal}. In
Section~\ref{se:final} we make some remarks that put our result in context.\bigskip

\noindent
{\bf Acknowledgement.} The first author thanks the University of Warwick
for its hospitality during January 2022, while this paper was in
progress. The research of the second author and the visit of the first
author were partly funded by
EPSRC grant EP/P031080/1.

\section{Characteristic zero}\label{se:0}

We include the well known proof for the connected case 
in characteristic zero for the sake
of convenience of comparison. It uses the
fact that the algebra of rational cochains on a simply connected space
has a commutative model as a DG algebra.

\begin{theorem}\label{th:rational}
Let $G$ be a connected Lie group\index{Lie group} 
(or more generally, any path connected 
topological group\index{topological group} with finite
dimensional rational homology)
and let $\kk$ be a field containing $\bQ$, the field of
rational numbers. Then $C^*BG$ is formal.
\end{theorem}
\begin{proof}
In rational homotopy theory,\index{rational homotopy theory} 
if the Sullivan minimal model\index{Sullivan minimal model} of a space has zero
differential then the rational homotopy type is formal.
This is true for simply connected spaces whose rational cohomology is a
polynomial ring on even degree generators tensored with an exterior
algebra on odd degree generators. The reason is that we can choose 
arbitrary cocycles representing the generators of cohomology, and there
are no relations to satisfy apart from (graded) commutativity, which is automatic 
because of the commutativity of the minimal model.
In the case of a
connected Lie group, the rational cohomology is isomorphic
to the invariant ring\index{invariant ring} 
$H^*(BT;\bQ)^W$, where $T$ is a maximal (compact) torus\index{maximal torus} 
and $W$ is the Weyl group. This is a polynomial ring on
even degree generators. 
So $C^*(BG;\bQ)$ is formal, and hence so is $C^*(BG;\kk)$ by extension
of scalars.
For more details and background on rational homotopy theory, 
we refer the reader to Proposition~15.5 of 
F\'elix, Halperin and Thomas~\cite{Felix/Halperin/Thomas:2001a}, 
or Example~2.67 in the 
book of F\'elix, Oprea and Tanr\'e~\cite{Felix/Oprea/Tanre:2008a}. 
\end{proof}

\begin{remark}
In the theorem, for $G$ connected, $\Omega BG\Qhat\simeq G\Qhat$ has
homology an exterior algebra on odd degree generators, so it is again
formal. In fact, the proof shows that $BG\Qhat$ and $G\Qhat$ are both
intrinsically formal.\index{intrinsically formal}
\end{remark}

The obstruction to extending the proof to other characteristics is
that the algebra of mod $p$ cochains on a space usually does not have a commutative DG
algebra model. However,
we can give a modified version of this theorem using the 
method of internal gradings. We do this in several steps, in the rest of this paper.

\section{Finite approximations}\label{se:finite}

In this section, we show how to approximate compact Lie groups by 
finite groups. This is closely related to the work of 
Dwyer and Wilkerson~\cite{Dwyer/Wilkerson:1994a}, but better suited
to the internal grading method.

We begin with case of a circle group $T=S^1$. We have $H^*(BT;\bF_p)=\bF_p[x]$ 
with $|x|=-2$.
 Let $\mu_m$ be the 
finite subgroup of $T$ consisting of elements whose $m$th power is the
identity. Then assuming that $p^n\ge 3$, we have 
$H^*(B\mu_{p^n};\bF_p)=\bF_p[x_n] \otimes \Lambda(t_n)$,
with $|x_n|=-2$ and $|t|=-1$, and $x_n$ is the restriction of $x$.
The restriction map from $H^*(B\mu_{p^{n+1}};\bF_p)$ 
to $H^*(B\mu_{p^n};\bF_p)$ sends the element $x_{n+1}$ to $x_n$
and $t_{n+1}$ to zero. Let $\mu_{p^\infty}$ be the union of 
the chain of subgroups 
$\mu_p \subseteq \mu_{p^2} \subseteq \cdots$ of $T$, which we regard as a discrete
group isomorphic to $\bZ/p^\infty$. Then the restriction map $H^*(BT;\bF_p)\to 
H^*(B\mu_{p^\infty};\bF_p)$ is an isomorphism.

For a torus $T=(S^1)^r$ of rank $r$ it works similarly. The inclusion
of $\mu_{p^\infty}^r$ in $T$ induces an isomorphism
\[ H^*(BT;\bF_p) \to H^*(B\mu_{p^\infty}^r;\bF_p). \]

\begin{lemma}
Let $T$ be a torus of rank $r$ acted on by a finite group $W$, let
$\alpha\in H^2(W,T)$ 
be an element of order $m\ge 1$, and let $\mu_m^r$ be the
finite subgroup of $T$ consisting of elements whose $m$th power is the identity.
Then $\alpha$ is in the image of 
$H^2(W,\mu_m^r)\to H^2(W,T)$.
\end{lemma}
\begin{proof}
We have a short exact sequence
\[ 1 \to \mu_m^r \to T \xrightarrow{m} T \to 1. \]
which gives an exact sequence
\[ \dots \to H^1(W,T) \to H^2(W,\mu_m^r) \to H^2(W,T) \xrightarrow{m}
  H^2(W,T) \to \cdots \]
If $\alpha\in H^2(W,T)$ has order $m$ then it is in the kernel of 
multiplication by $m$ on $H^2(W,T)$, and hence it is in the image of 
of $H^2(W,\mu_m^r)$.
\end{proof}

\begin{theorem}\label{th:approx}
Let $G$ be a compact Lie group with normal maximal torus $T$
of rank $r$ and let $W=G/T$. Then the Bousfield--Kan
$p$-completion $BG\phat$ is homotopy
equivalent to $B\breve G\phat$, where $\breve G$ is a locally
finite group sitting in a short exact sequence
\[ 1 \to \breve T \to \breve G \to W \to 1 \]
with $\breve T$ isomorphic to $\mu_{p^\infty}^r$.
If $|W|$ is coprime to $p$ then this sequence splits.
\end{theorem}
\begin{proof}
The group $G$ sits in a short
exact sequence 
\[ 1 \to T \to G \to W \to 1. \]
This extension defines an element $\alpha$ of $H^2(W,T)$. Since $|W|$ annihilates
$H^2(W,T)$, the order $m$ of $\alpha$ is a divisor of $|W|$. So by the lemma,
$\alpha$ is in the image of $H^2(W,\mu_m^r) \to H^2(W,T)$. It follows that
if we quotient out by $\mu_m^r$, we have a split short exact sequence
\[ 1 \to T/\mu_m^r \to G/\mu_m^r \to W \to 1. \]
Choosing a splitting, and taking the inverse image in $G$, 
we obtain a subgroup $\tilde W$ of $G$
that sits in a diagram
\[ \xymatrix{1 \ar[r]& \mu_m^r \ar[r]\ar[d] & 
\tilde W \ar[r]\ar[d]& W \ar[r]\ar@{=}[d] & 1\\
1 \ar[r] & T \ar[r] & G \ar[r] & W \ar[r] & 1.} \]
Let $\tilde G$ be the subgroup of $G$ generated by $\tilde W$ 
and $\mu_{p^\infty}^r$, regarded as a discrete locally finite group,
and let $\tilde T$ be the discrete locally finite 
subgroup of $\tilde G$ generated by 
$\mu_m^r$ and $\mu_{p^\infty}^r$.
Then we have a diagram
\[ \xymatrix{1 \ar[r] & \tilde  T \ar[r] \ar[d] &
\tilde G \ar[r]\ar[d] & W \ar[r]\ar@{=}[d] & 1 \\
1\ar[r] & T \ar[r] & G \ar[r] & W \ar[r] & 1.} \]
The comparison map of spectral sequences
\[ \xymatrix{H^*(BW; H^*(BT;\bF_p)) \ar@{=>}[r]\ar[d]& 
H^*(BG;\bF_p)\ar[d] \\
H^*(BW;H^*(B\tilde T;\bF_p)) \ar@{=>}[r] & H^*(B\tilde G;\bF_p)} \]
is an isomorphism on the $E_2$ page, and hence the map 
$\tilde G \to G$ induces
an isomorphism $H^*(BG;\bF_p) \to H^*(B\tilde G;\bF_p)$. 
It therefore induces a homology equivalence 
$H_*(B\tilde G;\bF_p) \to H_*(BG;\bF_p)$ and
a homotopy equivalence of Bousfield--Kan completions
$B\tilde G\phat \to BG\phat$.

Now the $p'$ elements of $\tilde T$ are the $p'$ elements of
$\mu_m^r$, and they form a characteristic 
finite subgroup $O_{p'}\mu_m^r$ of $\tilde T$, and hence also of
$\tilde G$.
Set $\breve G=\tilde G/O_{p'}\mu_m^r$ and
$\breve T=\tilde T/O_{p'}\mu_m^r\cong \mu_{p^\infty}^r$.
Then the quotient map $\tilde G \to \breve G$ is mod $p$
cohomology isomorphism, and $\breve G$ sits in a short
exact sequence
\[ 1 \to \breve T \to\breve G \to W \to 1. \]
The space $B\breve G\phat$ is homotopy equivalent to $BG\phat$.

For the final statement about splitting, the group $H^2(W,\mu_{p^n}^r)$
is annihilated by by $|W|$ and by $p^n$, and is hence zero. So
$H^2(W,\breve T)\cong H^2(W,\mu_{p^\infty}^r)\cong
\displaystyle\lim_{\substack{\longrightarrow\\n}}H^2(W,\mu_{p^n}^r)=0$.
\end{proof}

We shall make use of the theorem by means of the following proposition.

\begin{proposition}
Let $G$ be a compact Lie group with maximal torus $T$, and let
$W=N_G(T)/T$. If $p=0$ or $(p,|W|)=1$ then the inclusion $N_G(T)\to G$
induces a homotopy equivalence of $p$-completions
$BG\phat\simeq BN_G(T)\phat$.
\end{proposition}
\begin{proof}
In Theorem~20.3 of Borel~\cite{Borel:1967a}, it is proved that
for $G$ connected, $H^*(BG;\bQ) \to H^*(BT;\bQ)^W$ is an isomorphism.
Theorem~I.5 of Feshbach~\cite{Feshbach:1981a}, improves this 
to the statement
that whether or not $G$ is connected, the inclusion $T \to G$ induces
an isomorphism
\[ H^*(BG;\bZ) \otimes \bZ[1/|W|] \to H^*(BT,\bZ)^W \otimes
  \bZ[1/|W|]. \]
If $(p,|W|)=1$, then it follows using the five lemma on the 
long exact sequence in cohomology associated to
the short exact sequence of coefficients
\[ 0 \to \bZ[1/|W|] \xrightarrow{p} \bZ[1/|W|] \to \bF_p \to 0 \]
that $H^*(BG;\bF_p) \to H^*(BT;\bF_p)^W$ is an
isomorphism. This applies just as well to $N_G(T)$, and so
$H^*(BG;\bF_p) \to H^*(BN_G(T);\bF_p)$ is an isomorphism. Now complete
at $p$.
\end{proof}

\section{Formality}\label{se:formal}

Let $G$ be a compact Lie group with maximal torus $T$ of rank $r$, and
let $W=N_G(T)/T$. According to Theorem~\ref{th:approx},
there is a locally finite discrete group $\breve G$ sitting in
a short exact sequence
\[ 1 \to \breve T \to \breve G \to W \to 1 \]
with $\breve T \cong \mu_{p^\infty}^r$, and a homotopy equivalence 
$B\breve G\phat\simeq BG\phat$.

For the proof of formality, we shall need to make use of K\H{o}nig's
lemma from graph theory, so we begin by stating this.

\begin{lemma}
Let $\Gamma$ be a locally finite, connected infinite graph. Then $\Gamma$
contains a ray, namely an infinite sequence of vertices without repetitions
$v_0,v_1,v_2,\dots$ such that there is an edge from each $v_i$ to $v_{i+1}$.
\end{lemma}
\begin{proof}
This is proved in K\H{o}nig~\cite{Konig:1927a}. 
For a modern reference in English, see for example
Lemma~8.1.2 of Diestel~\cite{Diestel:2017a}.
\end{proof}

\begin{theorem}\label{th:formal}
If $(p,|W|)=1$ then
the $A_\infty$ algebra $C^*(B\breve G;\bF_p)$ is formal.
\end{theorem}
\begin{proof}
We first examine the case $\breve G=\breve T=\mu_{p^\infty}$.
Regarding $S^1$ as the unit circle  
in $\bC$, let $g_n=e^{2\pi\bi/p^n}$ be a generator for $\mu_{p^n}$ 
as a subgroup of $S^1$, so that $g_n^p=g_{n-1}$. So
$\mu_{p^\infty}$ is the union of these subgroups, regarded as an  
infinite discrete group.  
 
Set $X_n=g_n-1$, a generator for the radical $J(\bF_p\mu_{p^n})$.
Then $\bF_p\mu_{p^n}=\bF_p[X_n]/(X_n^{p^n})$.
We put
an internal $\bZ[\frac{1}{p}]$-grading on the group algebra $\kk\mu_{p^n}$ by setting
$|X_n|=1/p^n$. Then the fact that $X_n^p=X_{n-1}$ implies that
the inclusion $\bF_p\mu_{p^n}\hra\bF_p\mu_{p^{n+1}}$ preserves the internal
grading. Thus the union $\bF_p\mu_{p^\infty}$ is
$\bZ[\frac{1}{p}]$-graded, with all degrees lying in the interval
$[0,1)$.

Assuming that $p^n\ge 3$, we 
have $H^*(B\mu_{p^n};\bF_p)=\bF_p[x_n]\otimes \Lambda(t_n)$. 
This ring
inherits an internal grading from the group algebra, and we have
$|x_n|=(-2,-1)$ and $|t_n|=(-1,-\frac{1}{p^n})$.
The restriction map
from $H^*(B\mu_{p^{n+1}};\bF_p)$ to $H^*(B\mu_{p^n};\bF_p)$ sends $x_{n+1}$ to $x_n$
and $t_{n+1}$ to zero. Therefore
$H^*(B\mu_{p^\infty};\bF_p) = \kk[x]$
with $|x|=(-2,-1)$.

Now the $A_\infty$ structure maps $m_i$ on $H^*(B\mu_{p^\infty};\bF_p)$ given by
Kadeishvili's Theorem~\cite{Kadeishvili:1982a} preserves the internal grading, and
increase the homological grading by $i-2$. But for every non-zero
element of $H^*B\bZ/p^\infty$, the homological grading is twice the
internal grading. So for $m_i$ to be non-zero, we must have
$i=2$. Thus $C^*(B\mu_{p^\infty};\bF_p)$ is a formal $A_\infty$ algebra.

Similarly, for a larger rank torus, $\breve T\cong\mu_{p^\infty}^r$, 
we put an internal $\bZ[\frac{1}{p}]$-grading on $\bF_p\mu_{p^\infty}^r$ by adding the
internal gradings on the factors $\mu_{p^\infty}$, so that all degrees lie in the
interval $[0,r)$. This puts an internal grading on $H^*(B\mu_{p^\infty}^r;\bF_p)$, so that it
is a polynomial ring on $r$ indeterminates, all in degree $(-2,-1)$.
The same argument as in the case of $\mu_{p^\infty}$ now implies that 
$C^*(B\breve T;\bF_p)$ is formal.

Next, we come to the case where $\breve T\cong\mu_{p^\infty}^r$ 
is normal in $\breve G$ and $W=\breve G/\breve T$ is a
finite $p'$-group. In this case, the extension
\[ 1\to \breve T \to \breve G \to W \to 1\]
splits by Theorem~\ref{th:approx}.

In this case, we need to be more careful about the invariance under $W$ of the
grading. The group $W$ acts on the subgroup $\mu_p^r$ of elements of
$T$ of order $p$, and this defines an $\bF_pW$-module $M$ of dimension
$r$. The action of $W$ on $J(\bF_p\mu_{p^n}^r)/J^2(\bF_p\mu_{p^n}^r)$ makes it
a $\bF_p W$-module canonically
isomorphic to $M$. Since $\bF_p W$ is semisimple,
the short exact sequence of $\bF_p W$-modules
\[ 0 \to J^2(\bF_p\mu_{p^n}^r) \to J(\bF_p\mu_{p^n}^r) \to
  J(\bF_p\mu_{p^n}^r)/J^2(\bF_p\mu_{p^n}^r) \to 0 \]
splits. 
Such a splitting gives us a $W$-invariant
grading on $\bF_p\mu_{p^n}^r$, by lifting a basis. 
The problem is that it is not unique,
so we need to worry about compatibility. 

Now the diagram
\[ \xymatrix{\bF_p\mu_{p^n}^r \ar@{^(->}[r]\ar[dr]_{x \mapsto x^p} & 
\bF_p\mu_{p^{n+1}}^r \ar@{->>}[d]\\ 
&\bF_p\mu_{p^n}^r} \]
commutes, because over $\bF_p$, we have 
$(\sum_i a_ix_i)^p=\sum_i a_ix_i^p$. So a splitting for $n+1$ gives a
splitting for $n$ by taking $p$th powers.
Since there are only finitely many
splittings at each stage, 
it follows using K\H{o}nig's lemma 
that we may choose consistent splittings for
all $n>0$. The graph to which the lemma is applied has as vertices the
pairs consisting of a value of $n$ and a splitting for the above
sequence. The edges go from the pairs with $n+1$ to the pairs with
$n$, by taking $p$th powers. A ray in this graph consists of a
consistent set of splittings for all $n>0$.

Let $X_{n,1},\dots,X_{n,r}$ be bases for such a consistent
set of splittings, so that $X_{n+1,i}^p=X_{n,i}$. We may now put a
grading on $\bF_p\mu_{p^\infty}^r$ in such a way that the degree of
$X_{n,i}\in \bF_p\mu_{p^n}^r$ is equal to $1/p^n$, and $W$ preserves
the grading. 

Now we can put a grading on $\bF_p\breve G$ by choosing a 
copy of $W$ in $\breve G$ complementary to $\breve T$, and
putting it in degree zero. Then 
\[ H^*(B\breve G;\bF_p) \cong H^*(B\breve T;\bF_p)^W \] 
inherits an internal grading, and the homological grading of
all elements is again twice the internal grading. So as before,
all $m_i$ with $i\ne 2$ in the $A_\infty$ structure on 
$H^*(B\breve G;\bF_p)$ are forced to be zero. It follows
that $C^*(B\breve G;\bF_p)$ is a formal $A_\infty$ algebra.
\end{proof}

\begin{proof}[Proof of Theorem~\ref{th:main}]
By Theorem~\ref{th:approx}, the $A_\infty$ algebra $C^*(BG;\bF_p)$
is quasi-isomorphic to $C^*(B\breve G;\bF_p)$,
which by Theorem~\ref{th:formal} is formal.
\end{proof}

\section{Final remarks}\label{se:final}

It follows from Stasheff and Halperin~\cite[Theorem~9]{Halperin/Stasheff:1970a} 
that if $X$ is a space and $R$ is a commutative ring of coefficients 
such that $H^*(X;R)$ is a polynomial algebra then $C^*(X;R)$ is formal.
For example, $H^*(BSU(n);\bZ)$ and $H^*(BU(n);\bZ)$ are polynomial
rings on Chern classes, and $H^*(BSp(n);\bZ)$ is a 
polynomial ring on Pontryagin classes,  so the algebras of cochains are formal for 
all commutative coefficients in these cases. Similarly, 
$H^*(BO(n);\bF_2)$ and
$H^*(BSO(n);\bF_2)$ are polynomial rings on Stiefel--Whitney classes,
so the algebras of cochains are formal over any commutative ring in which $2=0$.

In the case of a connected, simply connected compact Lie groups $G$,
Borel~\cite{Borel:1953a,Borel:1954a,Borel:1960a,Borel:1960b,Borel:1967a}
and Steinberg~\cite{Steinberg:1975a} have 
investigated torsion in $H^*(BG;\bZ)$. Borel comments at the end of 
Volume~II of his \emph{\OE uvres} that the
upshot of these papers is that
the following are equivalent:
\begin{enumerate}
\item $H^*(G;\bZ)$ has no $p$-torsion,
\item $H^*(G;\bF_p)$ is an exterior algebra on odd degree classes,
\item $H^*(BG;\bF_p)$ is a polynomial algebra on even classes,
\item $H^*(BG;\bZ)$ has no $p$-torsion,
\item $H^*(G/T;\bF_p)$ is generated by elements of degree two, where
  $T$ is a maximal torus,
\item Every elementary abelian $p$-subgroup is contained in a torus,
\item Every elementary abelian $p$-subgroup of rank at most three is contained in a torus,
\item The multiplicity of some fundamental root in the dominant coroot is divisible by $p$.
\end{enumerate}

The primes $p$ for which this occurs are called the 
\emph{torsion primes} for $G$. They are a subset of the 
\emph{bad primes}, which are those for which the multiplicity of some
fundamental root in the dominant root is divisible by $p$. But these
are not the same set. For example, if $G=Sp(n)$ then there are no
torsion primes, but $2$ is a bad prime. If $G=G_2$ then the
only torsion prime is $2$ but the bad primes are $2$ and $3$.

Putting these together, if $G$ is a connected, simply connected
compact Lie group, $\kk$ is a field of characteristic $p$, 
and $p$ is not a torsion prime then $H^*(BG;\kk)$ is a polynomial
ring on even classes, and $C^*(BG;\kk)$ is formal. 

Here is a table of the torsion primes and bad primes.
\begin{center}
\begin{tabular}{|ccc||ccc|} \hline 
type & bad & torsion & type & bad & torsion \\ \hline
$A_n$&$\varnothing$&$\varnothing$&$E_6$&$\{2,3\}$&$\{2,3\}$ \\
$B_n$&$\{2\}$&$\{2\}$&$E_7$&$\{2,3\}$&$\{2,3\}$ \\
$C_n$&$\{2\}$&$\varnothing$&$E_8$&$\{2,3,5\}$&$\{2,3,5\}$ \\
$D_n$&$\{2\}$&$\{2\}$&$F_4$&$\{2,3\}$&$\{2,3\}$\\
&&&$G_2$&$\{2,3\}$&$\{2\}$ \\ \hline 
\end{tabular}\medskip
\end{center}

If $G$ is a compact Lie group which is connected but not simply connected,
then there exists a torus $T^n$ and a simply connected group $H$, 
such that if $\pi$ is the torsion
subgroup of $\pi_1(G)$ then $G$ sits in a
short exact sequence
\[ 1 \to \pi \to T^n \times H \to G \to 1. \]
So $H^*(BG;\bZ)$ has $p$-torsion if and only if $H^*(G;\bZ)$ has
$p$-torsion, which happens if and only if either $p$ divides
$|\pi|$ or $p$ is a torsion prime for $H$.
Otherwise, $H^*(BG;\bF_p)$ is a polynomial ring in even
degree generators, and $C^*(BG;\bF_p)$ is formal. 
The case where $p$ divides $|\pi|$ is studied in
Borel~\cite{Borel:1954a},
Mimura and Toda~\cite{Mimura/Toda:1991a}.

For non-connected compact Lie groups, the situation is quite
different. For example, let $G$ be a semidirect product $T^2\rtimes
\bZ/2$, where $\bZ/2$ acts on the $2$-torus $T^2$ by inverting every element. Then
$H^*(BG;\bF_3)$ is isomorphic to $H^*(BT^2;\bF_3)^{\bZ/2}$. We have
$H^*(BT^2;\bF_3)=\bF_3[x,y]$ with $|x|=|y|=2$, and $H^*(BG;\bF_3)$ is the
subring generated by $x^2,xy,y^2$, which is not a polynomial ring.
On the other hand, $3$ is not a torsion prime for $H^*(BG;\bZ)$, and 
$C^*(BG;\bF_3)$ is formal by Theorem~\ref{th:main}.

For a finite group $G$ and a prime $p$ dividing $|G|$,  it follows from Benson and
Carlson~\cite{Benson/Carlson:1994a} that  
$H^*(BG;\bF_p)$ is a polynomial ring if and
only if $p=2$, $G/O(G)$ is an elementary abelian $2$-group, and
the polynomial generators are in degree one; 
$C^*(BG;\bF_2)$ is formal in this case. More generally, if $p=2$ and
$G$ has elementary abelian Sylow $2$-subgroups, then $C^*(BG;\bF_2)$
is formal, even though the cohomology ring does not have to be
polynomial. There are also other
formal examples. For example, it 
is shown in Benson~\cite{Benson:tame} that
$C^*(BM_{11};\bF_2)$ is formal, whereas $H^*(BM_{11};\bF_2)\cong
\bF_2[x,y,z]/(x^2y+z^2)$ with $|x|=3$, $|y|=4$, $|z|=5$ is not a
polynomial ring. It would be interesting to know whether there are any
formal examples in odd characteristic.

\bibliographystyle{amsplain}
\bibliography{../repcoh}

\end{document}